\documentclass[12pt]{amsart}

\usepackage{graphicx}
\usepackage{stmaryrd}
\usepackage{amssymb}

\usepackage{amsmath,amscd}

\usepackage[T1]{fontenc}
\usepackage{yfonts}

\newcommand{\M}{\mathcal{M}}
\newcommand{\N}{\mathcal{N}}

\renewcommand{\P}{\mathbb{P}}

\newcommand{\Bor}{\mbox{{Bor}}}

\renewcommand{\int}{\mbox{int}}

\theoremstyle{definition}

\newtheorem*{proposition}{Proposition}
\newtheorem*{question}{Question}

\author{Marcin Sabok} 

\address{Mathematical Institute,
  Wroc\l aw University, pl. Grunwaldzki $2\slash 4$,
  $50$-$384$ Wroc\l aw, Poland }

\email{sabok@math.uni.wroc.pl}

\title{A note on ccc forcings}

\begin{document}

\begin{abstract}
  The aim of this short note is to communicate a simple
  solution to the problem posed in \cite{Zpl:FI} as Question
  7.2.7: is it true that for every ccc $\sigma$-ideal $I$
  any $I$-positive Borel set contains modulo $I$ an
  $I$-positive closed set?
\end{abstract}

\maketitle

\section{Introduction}

We work in a fixed Polish space $X$ with a continuous,
strictly positive Borel measure $\mu$. We write $\Bor(X)$
for the family of Borel sets in $X$.

A $\sigma$-ideal in $X$ is a family $I\subseteq\Bor(X)$
which is closed under taking countable unions and subsets.
An associated Boolean algebra $\P_I$ is the quotient algebra
$\Bor(X)\slash I$. We say that a $\sigma$-ideal $I$ is ccc
if any antichain in $\P_I$ is countable.

The most common ccc $\sigma$-ideals are: the family of Borel
sets of measure zero, denoted by $\N$ and the family of
Borel meager sets, denoted by $\M$.

We say that a Borel set $B$ is $I$-positive if $B\not\in I$.
It is a well known fact that any $\N$-positive Borel set
contains a $\N$-positive closed set. On the other hand, for
any $\M$-positive Borel set $B$ there is a closed set $C$
for which $C\setminus B\in\M$, in other words $B$ contains
$C$ modulo $\M$. In \cite{Zpl:FI} the author asks if this
can be generalized to all ccc $\sigma$-ideals.

\begin{question}[\cite{Zpl:FI}, Question 7.2.7]
  Suppose that $I$ is a $\sigma$-ideal such that $\P_I$ is
  ccc. Is it true that every positive Borel set contains a
  positive closed set modulo the ideal $I$?
\end{question}

\section{The solution}

We will answer the above question negatively by giving an
example of a ccc $\sigma$-ideal $I$ and a $I$-positive Borel
set which does not contain modulo $I$ any closed
$I$-positive set.

Take $J=\M\cap\N$. It is clear that $J$ is ccc, so we only
need to prove the following proposition.

\begin{proposition}
  There exists an $J$-positive Borel set $A$ which does not
  contain modulo $J$ any $J$-positive closed set.
\end{proposition}
\begin{proof}
  Let $A$ be any Borel set such that $A\in\N\setminus\M$.
  Such a set can be obtained by decomposing $X$ into two
  sets $A\in\N$ and $B\in\M$. In particular $A$ is
  $J$-positive. Take any closed set $C$ such that
  $C\setminus A\in J$. Note that this implies that $C\in\N$.
  But also $C\in\M$ because any set of measure zero must
  have empty interior. This shows that $C\in J$.
\end{proof}

\end{document}